\documentclass[11pt]{amsart}

\numberwithin{equation}{section} \oddsidemargin=-.0cm
\usepackage{amssymb}
\usepackage{amsmath}
\usepackage{color}
\usepackage{graphicx, graphics, subfigure}
\usepackage{amsfonts, amssymb, amsmath, amsthm}

\usepackage{sidecap}
\sidecaptionvpos{figure}{c}

\usepackage[colorlinks=true, pdfstartview=FitV, linkcolor=blue,
            citecolor=blue, urlcolor=blue]{hyperref}  

\newtheorem{thm}{Theorem}[section]
\newtheorem{lem}{Lemma}[section]
\newtheorem{defi}{Definition}[section]

\newtheorem{prop}[thm]{Proposition}

\newtheorem{rem}{Remark}[section]
\newtheorem{cor}{Corollary}[section]

\def\bt{\begin{thm}}
\def\et{\end{thm}}
\def\bl{\begin{lem}}
\def\el{\end{lem}}
\def\bd{\begin{defi}}
\def\ed{\end{defi}}
\def\bc{\begin{cor}}
\def\ec{\end{cor}}
\def\bp{\begin{proof}}
\def\ep{\end{proof}}
\def\br{\begin{rem}}
\def\er{\end{rem}}

\def\bprop{\begin{prop}}
\def\eprop{\end{prop}}

\def\Forall{\text{ } \forall \:}
\def\d{\, \mathrm{d}}

\def\be{\begin{equation}}
\def\ee{\end{equation}}
\def\bes{\begin{equation*}}
\def\ees{\end{equation*}}
\def\bea{\begin{equation} \begin{aligned}}
\def\eea{\end{aligned} \end{equation}}
\def\beas{\begin{equation*} \begin{aligned}}
\def\eeas{\end{aligned} \end{equation*}}

\def\bi{\begin{itemize}}
\def\ei{\end{itemize}}
\def\ben{\begin{enumerate}}
\def\een{\end{enumerate}}

\def\xii{\xi}
\def\s{\mathfrak{s}}
\def \c{\mathfrak{c}}

\definecolor{rred}{rgb}{0.7,0,0.1}
\newcommand{\mkk}{\color{black}}
\newcommand{\mkr}{\color{black}}
\newcommand{\mkrr}{\color{black}}

\title[Non-Markovian Reduced Systems for SPDEs Driven by Additive Noise]{Non-Markovian Reduced Systems for Stochastic Partial Differential Equations: The Additive Noise Case}

\author[Micka\"el Chekroun]{Micka\"el D. Chekroun}
\address[MC]{Department of Mathematics, University of Hawaii at Manoa, Honolulu, HI 96822, USA, and 
Department of Atmospheric \& Oceanic Sciences, University of California, Los Angeles, CA 90095-1565, USA} 
\email{mdchekroun@math.hawaii.edu}
\email{mchekroun@atmos.ucla.edu}
\author[Honghu Liu]{Honghu Liu}
\address[HL]{Department of Atmospheric \& Oceanic Sciences, University of California, Los Angeles, CA 90095-1565, USA}
\email{hliu@atmos.ucla.edu}

\author[Shouhong Wang]{Shouhong Wang}
\address[SW]{Department of Mathematics, Indiana University, Bloomington, IN 47405, USA}
\email{showang@indiana.edu, \url{http://www.indiana.edu/~fluid}}



\begin{document}

\begin{abstract}

This article proposes for stochastic
partial differential equations (SPDEs) driven by additive noise, a novel approach for  the approximate parameterizations of the ``small" scales by the ``large" ones,  along with the  derivaton  of the corresponding reduced systems. 
This is accomplished by seeking for stochastic parameterizing manifolds (PMs) introduced in \cite{CLW13a} which are random manifolds aiming to provide \textemdash\,  in a mean square sense  \textemdash\, such approximate parameterizations. Backward-forward systems are designed to give access to such PMs as pullback limits depending through the nonlinear terms  on the time-history of the dynamics of the low modes when  the latter is simply approximated by its {\mkr stochastic linear component}. It is shown that the corresponding pullback limits can be efficiently determined, leading in turn to an operational procedure for the  derivation of non-Markovian reduced systems able to achieve good modeling performances in practice. This is illustrated on a stochastic Burgers-type equation, where it is shown that the corresponding non-Markovian features of these reduced systems play a key role to reach such performances.

\end{abstract}

\maketitle

\section{Introduction}

The reduction problem of stochastic partial differential equations (SPDEs) has attracted a lot of attention recently, and various approaches have been proposed, which include the amplitude equations approach \cite{BM13, Pradas_al12} and  the manifold-based approaches \cite{CLW13a, Kan_al12} among many others \cite{EMS01, FS09, GKS04, Rob08}; see also the references therein. 

In this article, we {\mkk extend to SPDEs driven by additive noise,} the strategy introduced in \cite{CLW13a} to derive effective non-Markovian  reduced equations. This approach is based on approximate parameterizations of the ``small" scales by the ``large" ones {\it via} the concept of stochastic parameterizing manifolds (PMs), where the latter are random manifolds aiming to improve in mean square error the partial knowledge of the full SPDE solution when compared with {\mkrr its projection} onto the resolved modes. 

{\it Backward-forward systems} are designed to give access to such PMs in practice. The key idea consists here of representing  the modes with high wave numbers  as a pullback limit depending on the time-history of the modes with low wave numbers. The resulting manifolds obtained by such a procedure are not subject to a spectral gap condition such as encountered in the classical stochastic invariant/inertial manifold theory. Instead, certain PMs can be determined under weaker non-resonance conditions; {\mkrr see \eqref{NR} below}. 

Such an idea {\mkk of parameterizing the high modes as a functional of the the time-history of the low modes}, has been used in the context of 2D-turbulence  \cite{EMS01}\footnote{{\mkrr We mention also data-based approaches  such as \cite{CKG11,kondrashov2013low} where it  has been  shown  \textemdash\, in the context of climate dynamics \textemdash\,  that an appropriate conditioning of the signal's high-frequency variability on the time-history of the low-frequency modes helps predict the path of the future variations.}}, with the essential difference  that   the pullback limits considered here are associated with backward-forward systems that are {\it partially coupled} in the sense that only the (past values) of the resolved variables force {\mkr nonlinearly} the equations of the unresolved variables, without any feedback in the dynamics of the former; see Eqns.~\eqref{LLL} below.  {\mkr  In such systems, the dynamics of the low modes  is simply approximated by its stochastic linear component which helps  simplify  its  integration when compared with the {\it fully coupled} versions\footnote{directly rooted in the work of \cite{DPD96}.} encountered previously for the approximation of stochastic inertial manifolds \cite{Kan_al12}, while still making possible  the achievement of a good parameterizing quality.}

Non-Markovian stochastic reduced systems are then derived based on such a PM approach. The reduced systems take the form of stochastic differential equations (SDEs) involving random coefficients that convey memory effects via the history of the Wiener process, and arise from the nonlinear interactions between the low modes, embedded in the ``noise bath."  These random coefficients exhibit typically an exponential decay of correlations whose rate depends explicitly on gaps arising in the non-resonances conditions. It is shown on a stochastic Burgers-type equation, that such PM-based reduced systems can achieve very good performance in reproducing the SPDE dynamics projected onto the resolved modes.

\section{Functional framework} \label{Sec_framework} 

Our functional framework takes place in a pair of Hilbert spaces  ($\mathcal{H}_1$,$\mathcal{H}$), such that  $\mathcal{H}_1$  is compactly and densely embedded in $\mathcal{H}$. Let $A:\mathcal{H}_1\rightarrow \mathcal{H}$ be a sectorial operator \cite[Def.~1.3.1]{Hen81} such that $-A$ is stable in the sense that its spectrum satisfies $\mathrm{Re} (\sigma(-A) )< 0$.  To deal with  nonlinear SPDEs for which the nonlinear terms are responsible of a loss of regularity compared to the ambient space $\mathcal{H}$, we consider  standard interpolated spaces $\mathcal{H}_\alpha$ between $\mathcal{H}_1$ and $\mathcal{H}$ (with $\alpha \in [0, 1))$\footnote{depending on the problem at hand; see {\it e.g.} \cite{Hen81}.} along with  perturbations of the linear operator $-A$ given by a one-parameter family, $\{B_{\lambda}\}_{\lambda \in \mathbb{R}}$, of bounded linear operators from  $\mathcal{H}_\alpha$ to  $\mathcal{H}$,  depending continuously on $\lambda$. Defining $L_{\lambda}$ as $-A+B_{\lambda},$ we are thus left 
with a one-parameter family of sectorial operators $\{-L_\lambda\}_{\lambda \in \mathbb{R}}$, each of them mapping  $\mathcal{H}_1$ into $\mathcal{H}$.

Our {\mkk main purpose} is to {\mkk present a general strategy\footnote{based on the notion of stochastic PMs introduced  in \cite{CLW13a} in the context of SPDEs driven by linear multiplicative noise.} able to provide effective reduced equations for nonlinear stochastic evolution equations driven by additive white noise:}
\begin{equation} \label{SEE}
\d u = (L_\lambda u + F(u)) \d t + \d W_t;
\end{equation}
where $F$ {\mkk is a continuous nonlinear map} from $\mathcal{H}_\alpha$ into $\mathcal{H}$ for some $\alpha \in [0,1)$. For simplicity, we {\mkk assume $L_\lambda$ to be self-adjoint with an orthonormal basis of eigenfunctions $\{e_k\}_{k\in \mathbb{N}}$ in  $\mathcal{H}$, with corresponding eigenvalues $\{\beta_k(\lambda)\}_{k\in \mathbb{N}}$. In our setting, only the first $N$ modes will be randomly forced according to
\be
W_t(\omega) = \sum_{i=1}^{N} \sigma_i  W_t^i(\omega)e_i, \qquad t \in \mathbb{R}, \; \sigma_i>0,\; \omega \in \Omega,
\ee
where $\{W^i\}_{1 \leq  i \leq N}$  is a family of  mutually independent standard two-sided Brownian motions, with paths in $C_0(\mathbb{R},\mathbb{R})$.  A natural space of realizations $\Omega=C_0(\mathbb{R},\mathbb{R}^N)$, can be then associated with $W$; $\Omega$ being endowed with its corresponding  Borel $\sigma$-algebra $\mathcal{F}$, its  filtration $\{\mathcal{F}_t\}$, and the Wiener measure $\mathbb{P}$; see \cite[Appendix A]{Arnold98}.
Throughout this article we will adopt the formalism of random dynamical systems (RDSs) \cite{Arnold98, CLW13a, CSG11}.}

We will assume hereafter that {\mkk the stochastic evolution equation \eqref{SEE}} is such that for any initial datum $u_0\in H_{\alpha}$ there exists a unique global pathwise solution of \eqref{SEE}; see {\it e.g.} \cite{DPZ08} for conditions.

\section{Stochastic Parameterizing Manifolds: Analytic Expressions and {\mkk ``Past-Noise'' Dependence}} \label{Sec_PM}

In this section,  following \cite{CLW13a},  stochastic parameterizing manifolds (PMs) are introduced. Backward-forward systems are designed to give access to such PMs in practice. The key idea consists of representing the modes with ``high'' wave numbers (as parameterized by the sought PM) as a pullback limit depending on the time-history of the modes with ``low'' wave numbers. The cut between what is ``low'' and what is ``high'' is organized as follows. The subspace $\mathcal{H}^{\c} \subset \mathcal{H}$ defined by,
\bea \label{Hc}
\mathcal{H}^{\c} := \mathrm{span}\{e_1, \cdots, e_m\}, 
\eea 
spanned by the $m$-leading modes will be considered as our subspace associated with the low modes. 
Its topological complements, $\mathcal{H}^{\s}$ and $\mathcal{H}^{\s}_\alpha$ in  respectively, $\mathcal{H}$ and $\mathcal{H}_\alpha$, will be considered as associated with the high modes. 
We will use $P_{\c}$ and $P_{\s}$ to denote the canonical projectors associated with $\mathcal{H}^{\c}$ and $\mathcal{H}^{\s}$.

\subsection{\bf Stochastic parameterizing manifolds: Definition} 

{\mkk For a given SPDE of type \eqref{SEE}, based on \cite[Sect.~8]{CLW13a},   we define a stochastic parameterizing manifold  $\mathfrak{M}$, as the graph of a random function $h^{\mathrm{pm}}$ from  $\mathcal{H}^{\c}$ to $\mathcal{H}^{\s}_\alpha$ which is aimed to provide, for any $u(t,\omega)$ solution of   \eqref{SEE} given for a realization $\omega$,  an approximate parameterization  of the ``high'' part $u_{\s}(t,\omega)=P_{\s} u(t,\omega)$ in terms of the  ``low'' one $u_{\c}(t,\omega)=P_{\c} u(t,\omega)$,  so that  the mean squared error, $\int_0^T \bigl \|u_{\s}(t,\omega) -  h^{\mathrm{pm}}(u_{\c}(t,\omega),\theta_t \omega) \bigr \|_\alpha^2 \, \d t $, is smaller than the variance of $u_{\s}$,  $\int_0^T  \|u_{\s}(t,\omega)\|_\alpha^2 \; \d t$, whenever $T$ is sufficiently large. In statistical terms, a PM function $h^{\mathrm{pm}}$ is thus such that the fraction of variance of $u_{\s}(t,\omega)$ unexplained by $h^{\mathrm{pm}}(u_{\c}(t,\omega),\theta_t \omega)$ is less than the unity; see \cite[Sect.~8]{CLW13a} for more details. More precisely we have}

\bd  \label{def:PM}

A stochastic  manifold $\mathfrak{M}$ {\mkk given by}
\beas
\mathfrak{M}(\omega) := \{\xii + h^{\mathrm{pm}}(\xii, \omega) \mid \xii \in \mathcal{H}^{\c}\},  \; \, \omega \in \Omega,
\eeas
with $h^{\mathrm{pm}}: \mathcal{H}^{\c} \times \Omega \rightarrow \mathcal{H}^{\s}_\alpha$ being a measurable mapping, is called a stochastic parameterizing manifold (PM) associated with the SPDE \eqref{SEE} (for some fixed $\lambda$)  if the following conditions are satisfied:

\bi
\item[(i)] The function $h(\cdot, \omega): \mathcal{H}^{\c} \rightarrow \mathcal{H}^{\s}_\alpha$ is continuous for each $\omega$. 
\item[(ii)] For any $u_0 \in \mathcal{H}_\alpha $, there exists a positive random variable $\omega \mapsto T_0(\omega; u_0)$, such that:
\bea \label{PM condition}
\int_0^T  \! \bigl \|u_{\s}(t,\omega;u_0) -  h^{\mathrm{pm}}(u_{\c}(t,\omega; u_0),\theta_t \omega) \bigr \|_\alpha^2 \, \d t   < \int_0^T \! \|u_{\s}(t,\omega;u_0)\|_\alpha^2 \, \d t, \Forall \omega \in \Omega, \; T > T_0(\omega; u_0),
\eea
where {\mkk $u_{\c}(t,\omega;u_0)=P_{\c} u(t,\omega;u_0)$ and $u_{\s}(t,\omega;u_0)=P_{\s} u(t,\omega;u_0)$, with  $u(t,\omega;u_0)$ being the solution of the SPDE ~\eqref{SEE} emanating from $u_0$, and driven by the realization $\omega$.}

\ei

For a given realization $\omega$ and a given initial datum $u_0$, if $u_{\s}(\cdot,\omega;u_0)$ is not identically zero, the parameterization defect of $\mathfrak{M}$  is defined as the following time-dependent ratio:
\be\label{Eq_QualPM}
Q(T, \omega; u_0) :=  \frac{\int_0^T  \bigl \|u_{\s}(t,\omega;u_0) - h^{\mathrm{pm}}(u_{\c}(t,\omega; u_0), \theta_t \omega) \bigr \|_\alpha^2 \, \d t }{ \int_0^T \|u_{\s}(t,\omega;u_0)\|_\alpha^2 \, \d t},  \quad T > T_0(\omega; u_0).
\ee

\ed

\subsection{\bf {\mkk Stochastic parameterizing manifolds  as pullback limits of} backward-forward systems}  {\mkk

 In this section, we consider the problem of determination of PMs for SPDEs such as given by Eq.~\eqref{SEE}.  Theorem \ref{THM_PM_analytic} below provides analytic solutions to this problem, when the nonlinearity $F$ is bilinear, denoted as $B$ below. 
 This is achieved by following the approach introduced in \cite{CLW13a} that we adapt to the case of additive noise. In that respect, we consider} the following {\it backward-forward system} associated with the SPDE \eqref{SEE}:
\begin{subequations}\label{LLL}
\begin{align}
& \mathrm{d} u^{(1)}_{\c} =  L_\lambda^{\c} u^{(1)}_{\c}\, \mathrm{d} s + \mathrm{d} P_{\c} W_s, 
 && s \in[ -T, 0],   \quad \; u^{(1)}_{\c}(s, \omega)\vert_{s=0} = \xii,  \label{LLL1} \\
& \mathrm{d} u^{(1)}_{\s} = \bigl( L_\lambda^{\s} u^{(1)}_{\s}  +  P_{\s} B(u^{(1)}_{\c}(s-T, \omega)) \bigr) \mathrm{d} s +  \mathrm{d} P_{\s} W_{s-T}, 
&&  s \in [0, T],  \qquad u^{(1)}_{\s}(s, \theta_{-T}\omega)\vert_{s=0}= 0, \label{LLL2}
\end{align}
\end{subequations}
where  $L_\lambda^{\c} := P_{\c} L_\lambda$, $L_\lambda^{\s} := P_{\s} L_\lambda$, and $\xi \in \mathcal{H}^{
\c}$.

In the system above, the  initial value of $u_{\c}^{(1)}$ is prescribed in fiber $\omega$, and the initial value of $u_{\s}^{(1)}$ in fiber $\theta_{-T}\omega$. The solution of this system is obtained by using a {\it backward-forward integration procedure} made possible due to the partial coupling  {\mkk between the equations constituting the system} \eqref{LLL} where $u^{(1)}_{\c}$ forces the evolution equation of $u^{(1)}_{\s}$ but not reciprocally. {\mkk Note that since, $u_{\c}^{(1)}$ emanating (backward) from $\xi$ in $\mathcal{H}^c$, forces the equation ruling the evolution of $u_{\s}^{(1)}$, the latter depends naturally on $\xi$ and we will emphasize thus this dependence as $u_{\s}^{(1)}[\xi]$ hereafter. }{\mkk Theorem \ref{THM_PM_analytic} below  identifies {\it non-resonance  conditions}  \eqref{NR} under which the {\it pullback limit} of  $u_{\s}^{(1)}[\xi]$ exists; these non-resonance conditions differing interestingly  from those identified for the multiplicative noise case in \cite{CLW13a}.  

As it will be illustrated in the next section in the case of a stochastic Burgers-type equation driven by additive noise, such a pullback limit may be used in practice to give efficiently  access to PMs,  for a broad class of regimes; see also \cite{CLW13a}. The following theorem provides an analytical description of such PMs which in particular emphasizes the dependence on {\it the past of the noise path} of these manifolds but of a different form than arising in the multiplicative noise case \cite{CLW13a}. As in the latter case though, such features result from the nonlinear self-interactions of the low modes, ``embedded in the noise bath,''\footnote{The nature of this ``noise bath'' being additive compared to \cite{CLW13a}.}   
and such as projected onto the high modes, {\it i.e.}~from  $P_{\s} B(u^{(1)}_{\c}(s-T, \omega))$ in \eqref{LLL}; see the proof sketched below. }

\bt\label{THM_PM_analytic}
{\mkk Let us consider the SPDE \eqref{SEE} as given within the functional setting of Section \ref{Sec_framework}, with here $F$ assumed to be a bilinear function $B$.}
Assume also that $\beta_n(\lambda) < 0$ for all $n >  m$. Let ${\mkk \mathcal{I} }:= \{1, \cdots, m\}$ {\mkk with $m=\textrm{dim}(\mathcal{H}^\c)$. We assume furthermore that the following {\it non-resonance conditions} hold}:
\begin{equation}  \label{NR} \tag{NR}
\begin{aligned} 
& \Forall \, (i_1, i_2) \in \mathcal{I}^2,  \ n >  m, \quad  \text{ if } \quad  \langle B(e_{i_1},  e_{i_2}), e_n \rangle \neq 0,   \quad  \text{ then } \quad  \beta_{i_1}(\lambda)  + \beta_{i_2}(\lambda) - \beta_n(\lambda) > 0, \\
& \Forall \, (i_1, i_2) \in \mathcal{I}^2,  \ n >  m, \quad  \text{ if } \quad  \langle B(e_{i_1},  e_{i_2}), e_n \rangle \neq 0   \text{ and } \sigma_{i_1} \neq 0, \quad  \text{ then } \quad   \beta_{i_2}(\lambda) - \beta_n(\lambda) > 0, \\ 
& \Forall \, (i_1, i_2) \in \mathcal{I}^2,  \ n >  m, \quad  \text{ if } \quad  \langle B(e_{i_1},  e_{i_2}), e_n \rangle \neq 0   \text{ and } \sigma_{i_2} \neq 0, \quad  \text{ then }  \quad  \beta_{i_1}(\lambda)  - \beta_n(\lambda) > 0.
\end{aligned}
\end{equation}

Then the pullback limit of the solution $u_{\s}^{(1)}[\xii](T, \theta_{-T}\omega; 0)$ to \eqref{LLL2} exists and is given by:
\bea  \label{Eq_PBA}   
h^{(1)}_\lambda(\xii, \omega) & = \lim_{T \rightarrow +\infty} u_{\s}^{(1)}[\xii](T, \theta_{-T}\omega; 0) \\
&  = \int_{-\infty}^{0}  e^{-\tau L_\lambda^{\s}} L_\lambda^{\s} P_{\s} W_\tau(\omega) \d \tau + 
 \int_{-\infty}^{0} e^{-\tau L_\lambda^{\s}} P_{\s} B(u_{\c}^{(1)}(\tau, \omega; \xi)) \d \tau, \quad \Forall \, \xii \in \mathcal{H}^{\c}, \; \omega \in \Omega,
\eea
where $u_{\c}^{(1)}(s, \omega; \xi) $ is the solution of \eqref{LLL1}:
\begingroup   
\abovedisplayskip=0.2ex 
\belowdisplayskip=0.1ex 
\be \label{eq:uc}
u_{\c}^{(1)}(s, \omega; \xi) = e^{s L_\lambda^{\c}} \xi + P_{\c}W_s(\omega) - \int_s^0 e^{(s - \tau) L_\lambda^{\c}}  L_\lambda^{\c} P_{\c} W_\tau(\omega) \d \tau.
\ee
\endgroup
Moreover, $h^{(1)}_\lambda$ has the following {\mkk analytic expression}:
\bea  \label{eq:h1_expansion}
h^{(1)}_\lambda(\xi, \omega) & = \sum_{n=m+1}^{N} Z_{n,\lambda}(\omega)  e_n  + \sum_{n=m+1}^{\infty} \sum_{i_1 = 1}^m  \sum_{i_2 = 1}^m \Bigl(
A_{\lambda}^{n, i_1, i_2}(\omega) + B_{\lambda}^{n, i_1, i_2}(\omega) \xi_{i_1}  \\
& \hspace{3em} + C_{\lambda}^{n, i_1, i_2}(\omega) \xi_{i_2} +  D_{\lambda}^{n, i_1, i_2} \xi_{i_1} 
\xi_{i_2} \Bigr) \langle B(e_{i_1}, e_{i_2}), e_n \rangle e_n,
\eea
where $\xi_i  = \langle \xi, e_i \rangle$, $i = 1, \cdots, m$, and 
\beas
& Z_{n,\lambda}(\omega) := \sigma_n \beta_n(\lambda) \int_{-\infty}^{0} e^{-\tau \beta_n(\lambda)} W_\tau^{n}(\omega) \d \tau, \;\; A_{\lambda}^{n, i_1, i_2}(\omega) = \sum_{j = 1}^4 M_{j,\lambda}^{n, i_1, i_2}(\omega),  \; \; B_{\lambda}^{n, i_1, i_2}(\omega) = \sum_{j = 5}^6 M_{j,\lambda}^{n, i_1, i_2}(\omega),  \\  
& C_{\lambda}^{n, i_1, i_2}(\omega) = \sum_{j = 7}^8  M_{j,\lambda}^{n, i_1, i_2}(\omega), \quad D_{\lambda}^{n, i_1, i_2} := \int_{-\infty}^0 e^{(\beta_{i_1}(\lambda) + \beta_{i_2}(\lambda) - \beta_{n}(\lambda))\tau} \d \tau =   \frac{1}{ \beta_{i_1}(\lambda) + \beta_{i_2}(\lambda) - \beta_{n}(\lambda)},
\eeas
with
\beas
M_{1,\lambda}^{n, i_1, i_2}(\omega) & :=  \sigma_{i_1} \sigma_{i_2} \int_{-\infty}^0 e^{-\beta_n(\lambda) \tau} W_\tau^{i_1}(\omega) W_\tau^{i_2}(\omega) \d \tau,  \\
M_{2,\lambda}^{n, i_1, i_2}(\omega) & :=  - \sigma_{i_1}\sigma_{i_2} \beta_{i_2}(\lambda) \int_{-\infty}^0 \Bigl ( e^{- \beta_n(\lambda) \tau} W_\tau^{i_1}(\omega)  \int_{\tau}^0 e^{(\tau - \tau') \beta_{i_2}(\lambda)} W_{\tau'}^{i_2}(\omega) \d \tau' \Bigr) \d \tau, \\
M_{3,\lambda}^{n, i_1, i_2}(\omega) &  :=  - \sigma_{i_1}\sigma_{i_2} \beta_{i_1}(\lambda) \int_{-\infty}^0 \Bigl ( e^{- \beta_n(\lambda) \tau} W_\tau^{i_2}(\omega)  \int_{\tau}^0 e^{(\tau - \tau') \beta_{i_1}(\lambda) } W_{\tau'}^{i_1}(\omega) \d \tau' \Bigr) \d \tau, \\ 
M_{4,\lambda}^{n, i_1, i_2}(\omega) & :=  \sigma_{i_1}\sigma_{i_2}\beta_{i_1}(\lambda) \beta_{i_2}(\lambda) \int_{-\infty}^0 \Bigl ( e^{- \beta_n(\lambda) \tau} \int_{\tau}^0 e^{(\tau - \tau') \beta_{i_1}(\lambda) } W_{\tau'}^{i_1}(\omega) \d \tau' \int_{\tau}^0 e^{(\tau - \tau') \beta_{i_2}(\lambda) } W_{\tau'}^{i_2}(\omega) \d \tau'  \Bigr) \d \tau, \\
M_{5,\lambda}^{n, i_1, i_2}(\omega) & :=  \sigma_{i_2} \int_{-\infty}^0 e^{( \beta_{i_1}(\lambda)- \beta_n(\lambda)) \tau } W_\tau^{i_2}(\omega) \d \tau, \\ 
M_{6,\lambda}^{n, i_1, i_2}(\omega) & :=  - \sigma_{i_2}\beta_{i_2}(\lambda) \int_{-\infty}^0 \Bigl ( e^{(\beta_{i_1}(\lambda)- \beta_n(\lambda) )\tau} \int_{\tau}^0 e^{(\tau - \tau') \beta_{i_2}(\lambda) } W_{\tau'}^{i_2}(\omega) \d \tau' \Bigr) \d \tau,\\
M_{7,\lambda}^{n, i_1, i_2}(\omega) & :=  \sigma_{i_1} \int_{-\infty}^0 e^{( \beta_{i_2}(\lambda)- \beta_n(\lambda)) \tau } W_\tau^{i_1}(\omega) \d \tau, \\
 M_{8,\lambda}^{n, i_1, i_2}(\omega) & :=  - \sigma_{i_1}\beta_{i_1}(\lambda) \int_{-\infty}^0 \Bigl ( e^{(\beta_{i_2}(\lambda)- \beta_n(\lambda) )\tau} \int_{\tau}^0 e^{(\tau - \tau') \beta_{i_1}(\lambda) } W_{\tau'}^{i_1}(\omega) \d \tau' \Bigr) \d \tau.
\eeas

\et

\bp
(Sketch) First note that the solution of \eqref{LLL} can be formally obtained by using the variation-of-constants formula followed by an integration by parts performed to the resulting stochastic convolution terms.
{\mkk By doing so}, the solution to \eqref{LLL2} is given by $u_{\s}^{(1)}[\xii](T, \theta_{-T}\omega; 0)  = - e^{TL_\lambda^{\s}}P_{\s} W_{-T}(\omega) + \int_{-T}^{0}  e^{-\tau L_\lambda^{\s}} L_\lambda^{\s} P_{\s} W_\tau(\omega) \d \tau +  \int_{-T}^{0} e^{-\tau L_\lambda^{\s}} P_{\s} B(u_{\c}^{(1)}(\tau, \omega; \xi)) \d \tau$, where $u_{\c}^{(1)}(\cdot, \omega; \xi)$ is the solution of \eqref{LLL1} which takes the form of \eqref{eq:uc}. 

{\mkk Since $\beta_{n}(\lambda)<0$ for $n > m$,} the term $e^{TL_\lambda^{\s}}P_{\s} W_{-T}(\omega)$ converges to zero as $T\rightarrow +\infty$. {\mkk This implies that the pullback limit} of $u_{\s}^{(1)}[\xii]$ {\mkk exists and takes} the form given in \eqref{Eq_PBA} provided that the two integrals involved {\mkk in the expression of} $u_{\s}^{(1)}[\xii]$ {\mkk converge as  $T\rightarrow +\infty$}. The existence of $\int_{-\infty}^{0}  e^{-\tau L_\lambda^{\s}} L_\lambda^{\s} P_{\s} W_\tau(\omega) \d \tau$ {\mkk is a consequence of $\beta_n(\lambda) < 0$ for $n> m$}. The existence of $\int_{-\infty}^{0} e^{-\tau L_\lambda^{\s}} P_{\s} B(u_{\c}^{(1)}(\tau, \omega; \xi)) \d \tau$ results from the non-resonance condition \eqref{NR}. {\mkk This can be seen by expanding the bilinear term $P_{\s} B(u_{\c}^{(1)}(\tau, \omega; \xi))$, using  the expression of $u_{\c}^{(1)}$ given in \eqref{eq:uc}, and the fact that $L_\lambda$ is self-adjoint. 
This leads  to {\mkr three types of terms factorized by $ \langle B(e_{i_1}, e_{i_2}), e_n \rangle e_n$ in \eqref{eq:h1_expansion}: the constant terms, the linear ones and  the quadratic ones,} where the coefficients $D_{\lambda}^{n, i_1, i_2}$ and $M_{j,\lambda}^{n, i_1, i_2}(\omega)$ therein,  are ensured to exist due to the \eqref{NR}-conditon. The condition $\beta_n(\lambda) < 0$ for $n> m$, ensuring the existence of the Ornstein-Uhlenbeck processes $Z_{n,\lambda}$, we conclude then to the existence of $h^{(1)}_\lambda$ along with its analytic expression given in \eqref{eq:h1_expansion}.
}

\ep

Note that under similar assumptions {\mkk to} \cite[Theorem 8.4]{CLW13a}, {\mkk it can be shown furthermore} that $h^{(1)}_\lambda$  {\mkk  such as provided by the Theorem above, constitutes a PM function for the given SPDE \eqref{SEE}.}

Note also that the random coefficients $Z_{n,\lambda}$ and $M_{j,\lambda}^{n, i_1, i_2}$ given above exhibit decay of {\mkk correlations} as {\mkk it can be checked by similar calculations performed for the proof} of \cite[Lemma ~9.1]{CLW13a}. In particular, $Z_{n,\lambda}$, $M_{5,\lambda}^{n, i_1, i_2}$ and $M_{7,\lambda}^{n, i_1, i_2}$ are standard {\mkk Ornstein-Uhlenbeck processes}, which exhibit exponential decay of correlations with rates {\mkk respectively} given by $|\beta_n(\lambda)|$, $\beta_{i_1}(\lambda) - \beta_n(\lambda)$ and {\mkk $\beta_{i_2}(\lambda) - \beta_n(\lambda)$}.   

As we will {\mkk illustrate in the next section, the decay of correlations of the aforementioned coefficients are responsible for bringing {\it extrinsic memory effects}\footnote{in the sense of \cite{Hai09, HO07}.} in the $h^{(1)}_\lambda$-based stochastic reduced systems of  Eq.~\eqref{SEE}; see \eqref{SDE_h1} below.  The memory effects will} turn out to {\mkk play an} important role in {\mkk the performance achieved by such systems for the modeling of the SPDE  dynamics projected onto the low modes}.  {\mkk Unlike with the multiplicative noise case \cite{CLW13a}, such memory effects come with constant and linear terms which are absent compared to the deterministic case, where the resulting $h^{(1)}_\lambda(\xi)$  is reduced to $\sum_{n=m+1}^{\infty} \sum_{i_1 = 1}^m  \sum_{i_2 = 1}^m D_{\lambda}^{n, i_1, i_2} \xi_{i_1}  \xi_{i_2} \big \langle B(e_{i_1}, e_{i_2}), e_n \rangle e_n.$}

\vspace{2ex}

\section{{\mkk PM-based Non-Markovian Reduced Systems: Application to a Stochastic Burgers Equation}}

In this section, we consider the following stochastic initial-boundary value problem on the interval $(0, l)$:
\be \label{eq:Burgers}
\mathrm{d} u = \big( \nu u_{xx}  + \lambda u  - \gamma u  u_x\big) \mathrm{d} t + \mathrm{d} W_t(x,\omega),
\ee
supplemented with the Dirichlet boundary conditions, $u(0,t; \omega) = u(l,t; \omega) = 0$, $t\geq 0$, and initial condition $u(x, 0; \omega) = u_0(x)$, $x\in (0,l)$,
where $\nu, \lambda$ and $\gamma$ are positive parameters, and $u_0$ is some appropriate initial datum. This problem can be cast into the abstract form \eqref{SEE} with $\mathcal{H}:=L^2(0,l)$,  
and $\mathcal{H}_1:=H^2(0,l)\cap H_0^1(0,l)$, see {\it e.g.}~\cite{DDT94}. Here, the noise  $W_t(x,\omega)$ is taken to be $\sum_{i = 1}^{N} \sigma_i e_i(x) W_t^i(\omega)$ as {\mkk above} with $e_i(x)$ {\mkk denoting} the eigenmodes of the linear part {\mkk $L_{\lambda} u=\nu u_{xx} + \lambda u$}.  {\mkk 
Such a stochastic model is inscribed in the long tradition of study of the Burgers turbulence subject to random forces; see among many others \cite{BFK00,Frisch95}.}

Note that the eigenvalues of $L_\lambda$ are given by $\beta_n(\lambda) := \lambda - \frac{\nu n^2\pi^2}{l^2}$,  $n \in \mathbb{N}$, and the corresponding eigenvectors are $e_n(x) := \sqrt{2/l}\sin(n\pi x/l)$, $x\in(0,l)$.  We consider below the case where the subspace $\mathcal{H}^{\c}$ is  spanned by the first two eigenmodes, {\it i.e.}  $\mathcal{H}^{\c} := \mathrm{span}\{e_1, e_2\}$.

By projecting \eqref{eq:Burgers} against the {\mkk low} modes $e_1$ and $e_2$, we have $\mathrm{d} u_{\c} = \big( L_\lambda^{\c} u_{\c} +P_{\c} B (u_{\c} + u_{\s})\big) \mathrm{d}t + \mathrm{d} P_{\c} W_t$, where as before $u_{\c} := P_{\c} u$ with $P_{\c}$ being the canonical projector associated with $\mathcal{H}^{\c}$. By {\mkk replacing  $u_{\s}(t,\omega) := P_{\s} u(t,\omega)$ with  the pullback limit $h^{(1)}_\lambda (\xi, \theta_t \omega)$ given by \eqref{eq:h1_expansion}, we obtain the following reduced equation} 
\be \label{SDE_h1}
\mathrm{d} \xii = \Big( L_\lambda^{\c} \xii +P_{\c} B \big(\xii + h^{(1)}_\lambda(\xii, \theta_t\omega)\big)\Big) \mathrm{d}t + \mathrm{d} P_{\c} W_t,
\ee
{\mkk aimed to provide an approximation of the SPDE dynamics projected onto the low modes.

Recalling that the random coefficients $Z_{n,\lambda}$ and $M_{j,\lambda}^{n, i_1, i_2}$ contained in the expansion of $h^{(1)}_\lambda$ exhibit decay of correlations, we can conclude that {\it extrinsic memory effects} in the sense of \cite{Hai09, HO07}  are thus conveyed by the drift part of \eqref{SDE_h1},  making such reduced systems  {\it non-Markovian}. }

The analytic form of $h^{(1)}_\lambda(\xi, \omega)$ can be obtained from \eqref{eq:h1_expansion} by noting that in this example, the nonlinear interactions $B_{i_1i_2}^n := \langle B(e_{i_1}, e_{i_2}), e_n \rangle$  are given by $ B_{i_1i_2}^n = -\gamma \langle e_{i_1} ( e_{i_2})_x, e_n \rangle$, which take the following form: 
\be \label{nonlinear_interaction}
 B_{i_1i_2}^n = - \frac{\gamma i_2 \pi}{\sqrt{2}l^{3/2}}, \text{ if }  n = i_1 + i_2 ; \quad  B_{i_1i_2}^n= - \frac{ \gamma i_2 \pi \mathrm{sgn}(i_1-i_2)}{\sqrt{2}l^{3/2}}, \text{ if } n = |i_1-i_2|; \quad  B_{i_1i_2}^n = 
0, \text{ otherwise}.
\ee
In particular, for the subspace $\mathcal{H}^{\c}$ chosen above,  we have  $\langle e_{i_1} ( e_{i_2})_x, e_n \rangle = 0$ for any $n \ge 5$ and $i_1,i_2 \in\{1,2\}$.

{\mkk From a practical viewpoint, it is however cumbersome to use directly the analytic formula of $h^{(1)}_\lambda$ to determine the vector field $P_{\c} B\big(\xii + h^{(1)}_\lambda(\xii, \theta_t\omega)\big)$ as $\xi$ varies in $\mathcal{H}^{\c}$. Instead, we approximate $h^{(1)}_\lambda$  ``on the fly''  --- along a trajectory $\xi(t,\omega)$ of interest --- substituting $h^{(1)}_\lambda (\xi(t,\omega),\theta_t \omega)$ by $u^{(1)}_{\s}[\xi(t,\omega)](t+T, \theta_{-T}\omega; 0)$ as obtained by integration of the backward-forward system \eqref{LLL}, for $T$ chosen sufficiently large; cf.~\cite[Sect.~11]{CLW13a}. This  is much more manageable and leads naturally to consider the following substitutive reduced equation,}
\begin{subequations} \label{Reduced_fly_sys}
\begin{align}
&\mathrm{d} \xii_t = \Big( L_\lambda^{\c} \xii_t +P_{\c} B\big(\xii_t + u^{(1)}_{\s}[\xi(t,\omega)](t+T, \theta_{-T}\omega; 0)\big) \Big) \mathrm{d}t + \mathrm{d} P_{\c} W_t,    \quad \xi(0, \omega) = \phi, \quad t > 0, \label{eq:fly_xi} \\
& \hspace{-3em}\text{where $\phi$ is appropriately chosen according to the SPDE initial datum, and $u^{(1)}_{\s}[\xi]$ is obtained via}  \nonumber \\ 
& \mathrm{d} u^{(1)}_{\c} =   L_\lambda^{\c} u^{(1)}_{\c}  \mathrm{d} s +\mathrm{d} P_{\c} W_s, \hspace{13em} u^{(1)}_{\c}(s, \omega)\vert_{s=t} = \xi(t, \omega), \quad s\in [t-T, t],  \label{layered_uc1}\\
& \mathrm{d} u^{(1)}_{\s} = \bigl( L_\lambda^{\s} u^{(1)}_{\s} + P_{\s} B(u^{(1)}_{\c}(s-T, \omega)) \bigr) \mathrm{d} s  + \mathrm{d} P_{\s} W_{s-T}, \hspace{1.2em} u^{(1)}_{\s}(s, \theta_{-T}\omega)\vert_{s=t}= 0, \quad \; s\in [t, T+t],  \label{layered_us}
\end{align}
\end{subequations}
{\mkk which written in coordinate form, gives (cf.~\cite[Section 11.2]{CLW13a})},
\begin{subequations} \label{Reduced_AE}
\begin{align}
& \d \xi_1 =\Bigl( \beta_1(\lambda) \xi_1 + \frac{\gamma\pi}{\sqrt{2}l^{3/2}} \Bigl (\xi_1\xi_2 + \xi_2 y_3^{(1)} + \sum_{j = 3}^N y_j^{(1)} y_{j+1}^{(1)} \Bigr ) \Bigr)\d t  + \sigma_1 \d W^1_t,  \hspace{2em} t > 0,  \label{SDE:xi_1} \\ 
& \d \xi_2 =\Bigl( \beta_2(\lambda) \xi_2 + \frac{\sqrt{2} \gamma \pi}{l^{3/2}} \Bigl ( - \frac{1}{2} (\xi_1)^2 + \xi_1 y_3^{(1)} + \xi_2 y_4^{(1)} + \sum_{j = 3}^{N} y_{j}^{(1)} y_{j+2}^{(1)} \Bigr) \Bigr) \d t  + \sigma_2\d W^2_t, \hspace{1.3em} t > 0, \label{SDE:xi_2}\\
& \text{with} \hspace{4em} \xi_1(0, \omega) = \langle \phi, e_1 \rangle, \qquad \xi_2 (0, \omega) = \langle \phi, e_2 \rangle,   \label{SDE:xi_ini}
\end{align}
\end{subequations}
where $\xi_t = \xi_1(t,\omega) e_1 + \xi_2(t,\omega)e_2$ and $y_i^{(1)}$, $i =3, \cdots, N$, are obtained {\it via},
\begin{subequations} \label{Reduced_AE2}
\begin{align}
& \d y^{(1)}_1 = \beta_1(\lambda) y^{(1)}_1 \d s + \sigma_1 \d W^1_s,   \hspace{15em}  s \in [t-T, t],  \label{SDE_uc1_y1}\\ 
& \d y_2^{(1)} = \beta_2(\lambda) y_2^{(1)}  \d s + \sigma_2 \d W^2_s,  \hspace{15em} s \in [t-T, t], \label{SDE_uc1_y2} \\
& \d y_3^{(1)} =  \Bigl(  \beta_3(\lambda) y_3^{(1)} - \frac{3\gamma \pi}{\sqrt{2}l^{3/2}} y_1^{(1)}(s - T, \omega) y_2^{(1)}(s-T, \omega) \Bigr) \d s  + \sigma_3 \d W^3_{s-T},  \hspace{1em} s \in [t, T+t],  \label{SDE_us_y1} \\ 
& \d y_4^{(1)} = \Bigl( \beta_4(\lambda) y_4^{(1)} - \frac{\sqrt{2}\gamma \pi}{l^{3/2}} [y_2^{(1)}(s-T, \omega)]^2  \Bigr) \d s + \sigma_4 \d W^4_{s-T}, \hspace{2em} s \in [t, T+t], \label{SDE_us_y2}\\ 
& \d y_j^{(1)} = \beta_j(\lambda) y_j^{(1)} \d s + \sigma_j \d W^j_{s-T}, \hspace{2em} s \in [t, T+t], \qquad j  = 5, \cdots, N, \label{SDE_us_yn}\\ 
& \text{with} \hspace{1em} y^{(1)}_1(s, \omega)\vert_{s=t} = \xi_{1}(t,\omega), \quad y_2^{(1)}(s, \omega)\vert_{s=t} = \xi_2(t,\omega),  \quad  y_i^{(1)}(s, \theta_{-T}\omega)\vert_{s=t} = 0, \quad i = 3, \cdots, N. \nonumber
\end{align}
\end{subequations}

{\bf \noindent Numerical results}. We assess below the performances achieved by the reduced system \eqref{Reduced_AE}-\eqref{Reduced_AE2}, in the modeling of the pathwise SPDE dynamics on the $\mathcal{H}^{\c}$-modes, on one hand; and of the full spatio-temporal field $u$ obtained by direct simulation of the SPDE \eqref{eq:Burgers}, on the other. {\mkr The importance of the memory terms conveyed by $h^{(1)}_\lambda$ is assessed by comparison with the performances achieved by an averaged version of \eqref{Reduced_AE}, which consists of replacing $h^{(1)}_\lambda$ in the non-Markovian reduced system \eqref{SDE_h1} by its expected value $\mathbb{E}(h^{(1)}_\lambda)$:} 
\be \label{SDE_averaged}
\mathrm{d} \xii = \Big( L_\lambda^{\c} \xii +P_{\c} B \big(\xii +\mathbb{E} \big(h^{(1)}_\lambda(\xii, \theta_t\omega)\big) \big)\Big) \mathrm{d}t + \mathrm{d} P_{\c} W_t.
\ee

It can be shown that the {\mkk concerned random coefficients provided by Theorem \ref{THM_PM_analytic}  come each with zero-expectation, resulting in an agreement of $\mathbb{E}(h^{(1)}_\lambda)$  with the pullback limit associated with the corresponding deterministic backward-forward system. The coordinate form of \eqref{SDE_averaged} can then be obtained easily; see \cite{CLW13a}.}

For all the simulations below, the noise is taken to be $W_t(x, \omega) = \sum_{i = 1}^{10} \sigma_i e_i(x) W_t^{i}(\omega)$, with $\sigma_1 = \cdots = \sigma_{10} =: \sigma>0$. The values of $\sigma$ and other system parameters are specified in the captions of the figures given below;  {\mkk the numerical schemes adopted here are  those used in \cite[Sect.~11]{CLW13a} adapted to the additive noise case.} In all the numerical experiments that follow, the time step has been taken to be $\delta t =0.01$, and {\mkrr the mesh size $\delta x$ has been specified as indicated in the Figures captions below.}

\begin{figure}[hbtp]
   \centering
   \includegraphics[height=7.5cm, width=16.4cm]{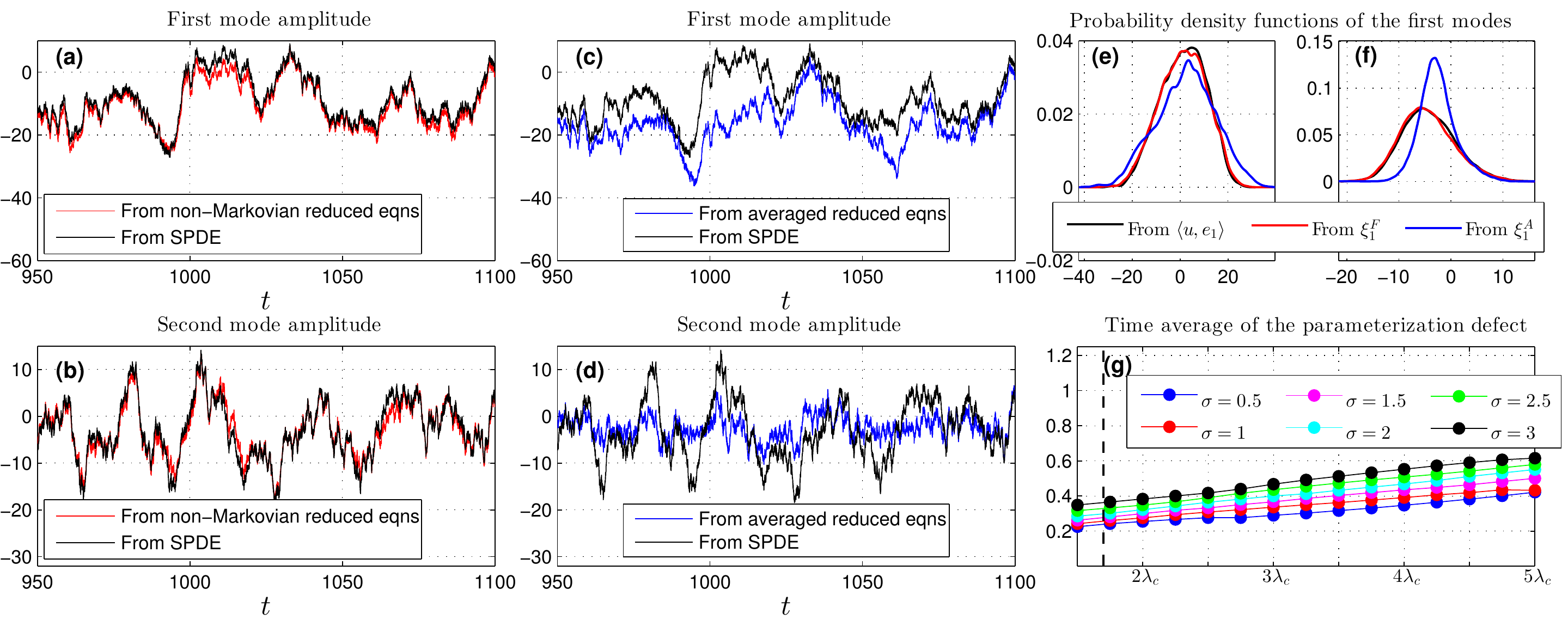}
  \caption[a]{{\footnotesize {\bf (a)-(b)}: The dynamics of the resolved modes modeled by the reduced system \eqref{Reduced_AE}-\eqref{Reduced_AE2} compared with the SPDE dynamics projected onto the resolved modes; {\bf (c)-(d)}: Analogous results given in {\bf (a)-(b)} produced by the averaged reduced system \eqref{SDE_averaged}; {\bf (e)-(f)}: Probability density functions (PDFs) of the first two modes amplitudes as estimated from a direct simulation of the SPDE \eqref{eq:Burgers} (black) {\it vs} those such as obtained by simulation of the reduced system \eqref{Reduced_AE}-\eqref{Reduced_AE2} (red), and its averaged   version \eqref{SDE_averaged} (blue); {\bf (g)}: Time averages of the parameterization defect of {\mkk $u^{(1)}_{\s}[\xi]$ for various $\lambda$ and $\sigma$}; here the time average {\mkrr (for a given realization $\omega$)} is $\frac{1}{T_2 - T_1}\int_{T_1}^{T_2} Q(T, \omega; u_0)\, \d T$, where $Q$ is computed using \eqref{Eq_QualPM} with $u^{(1)}_{\s}[u_{\c}(t,\omega)](t+T, \theta_{-T}\omega; 0)$ in place of $h^{\mathrm{pm}}_\lambda$; and $T_1$ and $T_2$ are taken here to be $400$ and $1000$, respectively. The vertical dashed line {\mkk on panel {\bf (g)}} corresponds to $\lambda = 1.7\lambda_c$ which is the value of $\lambda$ used for the results shown in panels {\bf (a)-(f)}. Here, $\lambda_c:= \nu \pi^2/l^2$ is the critical value at which the leading eigenvalue $\beta_1(\lambda)$ changes sign. {\mkk Except in {\bf (g)}, where the values of $\lambda$ and $\sigma$ are varied, the results reported here were computed for the following values of the parameters: $\gamma = 0.5$, $l = 7\pi$, $\nu = 2$, $\lambda = 1.7\lambda_c$, and $\sigma = 3$}. {\mkrr In all cases, $\delta x=0.1679.$} Similar to \cite{CLW13a}, it has been observed that the results reported here are robust w.r.t. the choice of the realization. }}   \label{fig1}
\end{figure}

Results of Fig.~\ref{fig1} are reported for a regime corresponding to a large amount of noise, and for a domain size $l$ corresponding to small decay {\mkrr rates} of correlations of the memory terms\footnote{For instance, the decaying rates of $Z_{3,\lambda}$ and $M_{5,\lambda}^{3,1,2}$ are given here by $|\beta_3(\lambda)| \approx 0.29$ and $\beta_1(\lambda) - \beta_3(\lambda) \approx 0.33$, respectively.}, so that in particular, such regimes correspond to a natural test bed for assessing the role of the extrinsic memory terms in the modeling performance achieved by the reduced system \eqref{Reduced_AE}-\eqref{Reduced_AE2}. This assessment is accomplished by comparing the results of \eqref{Reduced_AE}-\eqref{Reduced_AE2} with those of the averaged reduced system \eqref{SDE_averaged}. Figure \ref{fig1} (a)--(d) shows an episode where \eqref{Reduced_AE}-\eqref{Reduced_AE2} outperforms \eqref{SDE_averaged} in modeling the first two modes amplitudes. {\mkr It has been observed that such episodes repeat very often as time flows, 
as supported by Fig.~\ref{fig1} (e)-(f) from  the estimation of the corresponding probability density functions (PDFs). {\mkrr The latter show that the PDFs of the first two modes amplitudes as simulated by   \eqref{Reduced_AE}-\eqref{Reduced_AE2} coincide almost with those of the SPDE dynamics projected onto the $\mathcal{H}^{\c}$-modes, whereas the PDFs simulated from the averaged reduced system \eqref{SDE_averaged} fails in capturing the correct statistical behavior, demonstrating thus the modeling improvements brought by the non-Markovian features of \eqref{Reduced_AE}-\eqref{Reduced_AE2}.}}

\begin{figure}[!hbtp]
\centering
   \includegraphics[height=7cm, width=13cm]{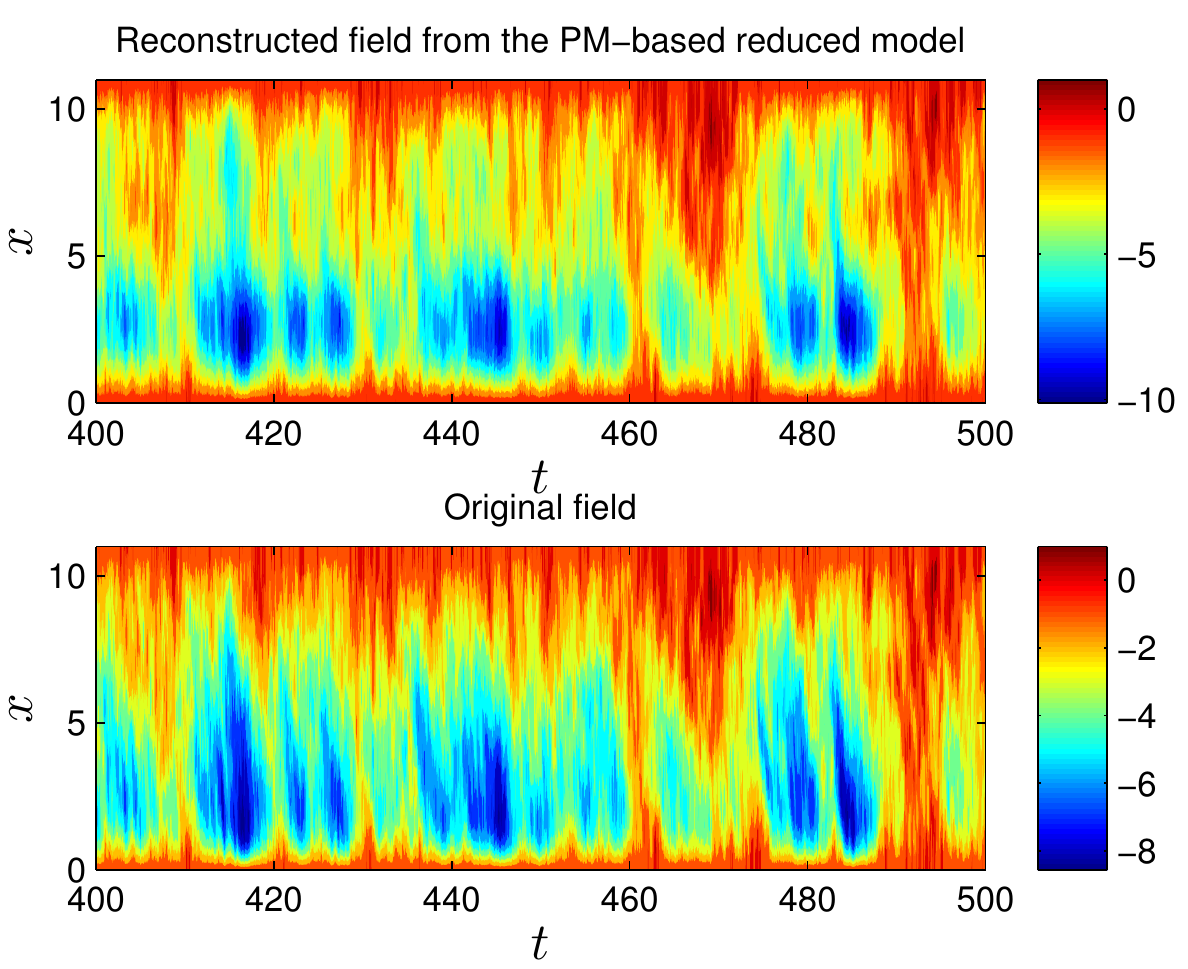}
  \caption{{\footnotesize {\mkk Spatio-temporal field, 
  $\xi(t,\omega)+u^{(1)}_{\s}[\xi(t,\omega)](t+T, \theta_{-T}\omega; 0)$, 
  reconstructed from the reduced system \eqref{Reduced_AE}-\eqref{Reduced_AE2} (upper panel).}
  {\mkk For comparison the SPDE solution field is plotted in the lower panel}. {\mkk Here the model parameters are $\gamma = 0.5$, $l = 3.5\pi$, $\nu = 2$, $\sigma = 1.5$, and $\lambda = 1.7\lambda_c$. {\mkrr Here $\delta x=0.0839.$}} 
  }}   \label{fig2}
\end{figure}

It has also been noticed that very good reconstructions of the SPDE solution $u$ can be obtained from the reduced model \eqref{Reduced_AE}-\eqref{Reduced_AE2}, in certain regimes.  Figure~\ref{fig2} illustrates such {\mkr a regime where the coarse-grained features of the spatio-temporal fluctuations of the solution $u$  are well reproduced from the solution $\xi(t,\omega)$ of the reduced model \eqref{Reduced_AE}-\eqref{Reduced_AE2} once the nonlinear corrective term $u^{(1)}_{\s}[u_{\c}(t,\omega)](t+T, \theta_{-T}\omega; 0)$ has been added.} 
For all the numerical {\mkk experiments}, the pullback time $T$ in \eqref{Reduced_AE}-\eqref{Reduced_AE2} is fixed to be $T=2$. It has been checked that {\mkk an increase of $T$ has little impact on the quality of the results reported here}.

\begin{figure}[hbtp]
   \centering
   \includegraphics[height=4cm, width=16cm]{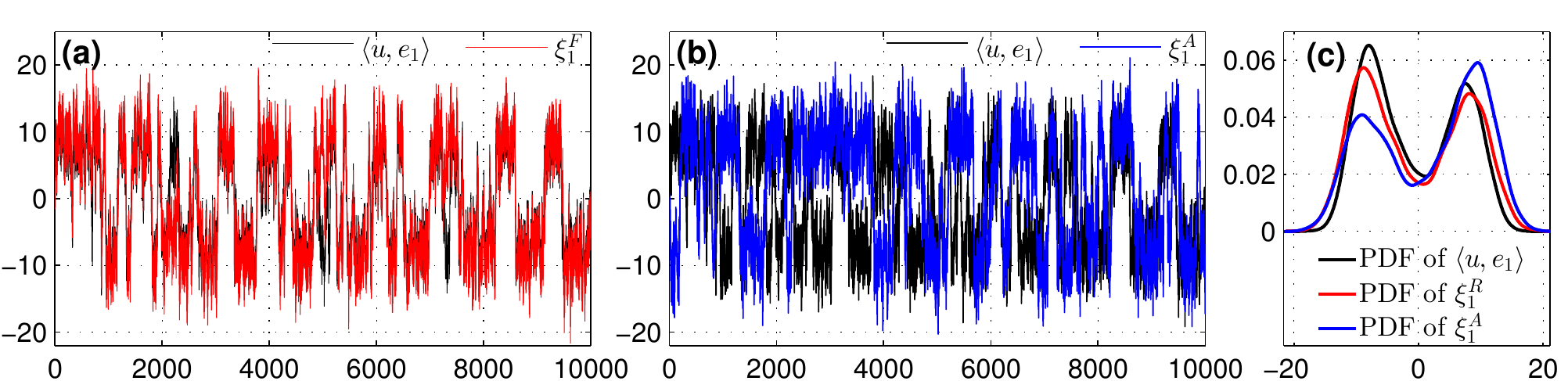}
  \caption{{\footnotesize {\bf (a)-(b) and (c)}: Bimodal behavior of the SPDE dynamics as projected onto $e_1$ along with its modeling from the non-Markovian system \eqref{Reduced_AE}-\eqref{Reduced_AE2} (red), and from its averaged version  \eqref{SDE_averaged} (blue). Noise-induced transitions occurring in time are well-captured by the non-Markovian reduced system when  compared to its averaged version.  Parameters are those of Fig.~\ref{fig2}.
}}   \label{fig3}
\end{figure}

Modeling performances such as reported in Fig.~\ref{fig1} (a)-(d) and Fig.~\ref{fig2}, are underpinned by the fact that $h^{(1)}_\lambda$ provides a good parameterization of the unresolved modes by the resolved ones as supported by {\mkk the computations of the (time-averaged) parameterization defects} shown in Fig.~\ref{fig1} (g). 
This is particularly remarkable given that a 
{\mkr pathwise SPDE solution $u$ may follow  a {\it bimodal behavior} such as illustrated in Fig.~\ref{fig3}  for the first mode amplitude.\footnote{for the same parameters values  used for Fig.~\ref{fig2}.}}

It has been observed that the transitions between the {\mkr corresponding statistical equilibrium states are very} often well captured by the reduced model \eqref{Reduced_AE}-\eqref{Reduced_AE2} with nevertheless some failures that may occur as time flows; see Fig.~\ref{fig3} (a). {\mkr The overall bimodal behavior as observed on the first mode is however well-captured by the non-Markovian system \eqref{Reduced_AE}-\eqref{Reduced_AE2} when compared to its averaged version  \eqref{SDE_averaged} (cf. Fig.~\ref{fig3} (c)), demonstrating once again the importance of the memory effects in reproducing non-trivial statistical features of the dynamics as well as the {\it noise-induced transitions} occurring in time; compare  Fig.~\ref{fig3} (a) with Fig.~\ref{fig3} (b). The ability of capturing such transitions accurately is tightly related to the modeling problem of the excitation of the small scales by the noise through the nonlinear terms.}

 The differences in the finer-scale details as observed in Fig.~\ref{fig2} are for instance due to such noise-induced phenomenon and the  limited number of modes resolved ({\mkk here only two}). In that respect, it is worthwhile to mention that $u^{(1)}_{\s}$  obtained by integration of \eqref{layered_us}, belongs to $\mathrm{span}\{e_3, \cdots, e_{10}\}$ with its components onto the fifth up to the tenth modes being modeled just as a red noise (see \eqref{SDE_us_yn}), due to the particular nonlinear self-interactions of the low modes for the problem at hand; see \eqref{nonlinear_interaction}. 
Such a ``standard'' red noise approximation for the fifth to the tenth modes is  {\mkk however}
not detrimental to the modeling performances achieved by the reduced system \eqref{Reduced_AE}-\eqref{Reduced_AE2}, when these modes do not contain a too important fraction, $\mathfrak{f}_{5}^{10}(T):=|\sum_{j=5}^{10} \langle u,e_j\rangle e_j|_{L^2(0,T;\mathcal{H}_{\alpha})}^2|u_{\s}|_{L^2(0,T;\mathcal{H}_{\alpha})}^{-2}$, of the unresolved variance of the signal, as $T$ evolves. {\mkk For such a case, the non-Markovian features brought by the more elaborated stochastic ``reddish'' processes
$M_{j,\lambda}^{n, i_1, i_2}$ ($n\in\{3,4\}$, $i_1, i_2 \in \{1,2\}$)  \textemdash\, arising typically with non-Gaussian statistics ({\it cf.}~\cite[Fig.~1]{CLW13a})  \textemdash\, in the reduced system \eqref{Reduced_AE}-\eqref{Reduced_AE2}, are sufficient to achieve a good parameterization of the excitation of the small scales by the noise,  leading to the successes reported in Figures \ref{fig1},   \ref{fig2} and \ref{fig3}.}

When $\mathfrak{f}_{5}^{10}$ happens to get large, the episodes of time where the modeling performances achieved by \eqref{Reduced_AE}-\eqref{Reduced_AE2} get deteriorated, become more numerous.  A remedy to such episodes of deterioration relies on the usage  of other parameterizing manifolds based on {\it multilayer backward-forward systems} such as introduced in \cite{CLW13a}. The latter convey a ``matriochka'' of nonlinear self-interactions of the low modes which arise typically with a {\it hierarchy of memory effects} of more elaborated structures than conveyed by $h^{(1)}_{\lambda}$; see \cite[Section 11.3]{CLW13a} for the multiplicative noise case. The corresponding results for the additive noise case will be reported elsewhere.

\section*{Acknowledgments}

MDC and HL are supported by the National Science Foundation
 grant DMS-1049114  and Office of Naval Research grant
N00014-12-1-0911. SW is supported in part by National Science Foundation grants DMS-1211218, and DMS-1049114, and by Office of Naval Research grant N00014-11-1-0404.

\bibliographystyle{amsalpha}
\newcommand{\etalchar}[1]{$^{#1}$}
\providecommand{\bysame}{\leavevmode\hbox to3em{\hrulefill}\thinspace}
\providecommand{\MR}{\relax\ifhmode\unskip\space\fi MR }
\providecommand{\MRhref}[2]{%
  \href{http://www.ams.org/mathscinet-getitem?mr=#1}{#2}
}
\providecommand{\href}[2]{#2}

\end{document}